\documentclass[english,a4paper,11pt]{article}
\usepackage{theorem}
\usepackage{amssymb,amsmath}
\usepackage{latexsym}
\usepackage{babel}
\usepackage{color}

\usepackage{color,graphicx,array}

\newcommand{\be}{\begin{equation}}
\newcommand{\ee}{\end{equation}}

\newtheorem{thm}{Theorem}[section]

\newtheorem{lem}{Lemma}[section]
\newtheorem{rem}{Remark}[section]

\newtheorem{defin}{Definition}[section]

\newtheorem{exam}{Example}[subsection]
\DeclareFontFamily{T1}{cmr}{\hyphenchar\font=-1}

\numberwithin{equation}{section}
\DeclareFontFamily{T1}{cmr}{\hyphenchar\font=-1}
\usepackage{color}
\usepackage{latexsym}
\usepackage[numbers,sort&compress,square]{natbib}

\begin{document}

 \noindent \textbf{{\LARGE Approximation of the fixed point of the product of two operators in Banach algebras with applications to some functional equations}}

 \noindent \textbf{Khaled Ben Amara$^{1}$, Maria Isabel Berenguer$^{2,*}$ Aref Jeribi$^{3}$}

 \noindent $^{(1,3)}$  Department of Mathematics. Faculty of Sciences of Sfax.
University of Sfax. Road Soukra Km $3.5 B.P. 1171, 3000,$ Sfax,
Tunisia. \\
 $^{(2)}$ Department of Applied Mathematics and Institute of Mathematics (IEMath-GR),  E.T.S. de Ingenieria de Edificaci\'on,  University of Granada, Granada, Spain.
\\
 $^{(*)}$ Corresponding author.

\date{}

\centerline {e-mails: $^{(1)}$ khaled.benamara.etud@fss.usf.tn,
$^{(2)}$ maribel$$@$$ugr.es, $^{(3)}$ aref.jeribi$$@$$fss.rnu.tn}

\renewcommand{\thefootnote}{}

\noindent {\bf\small Abstract.}
In this paper, the existence and uniqueness of the fixed point for the product of two nonlinear operator in Banach algebra is discussed. In addition, an approximation method of the fixed point of hybrid nonlinear equations  in Banach algebras is established. This method is applied to two interesting different types of functional equations. In addition, to illustrate  the applicability of our  results we give some numerical examples.

\vskip0.7cm
\section{\bf{Introduction}}
Many phenomena in physics,
chemistry, mechanics, electricity, and so as, can be formulated
by using the following nonlinear differential equations with a nonlocal initial condition of the form:

\begin{equation}\label{aplicationmulti1a}\left\{
  \begin{array}{rl}
 \displaystyle \frac{d}{dt}\left(\frac{x(t)}{f(t,x(t))}\right)&=  g(t,x( t)), t\in J,\\\\
   x(0)&= \Gamma(x),
 \end{array} \tag{P1}
\right.\end{equation}
where $\Gamma$ is a mapping from $C(J)$ into $\mathbb{R}$ which represents the nonlocal initial condition of the considered problem, see \textup{\cite{Deimling,Djebali}}. The nonlocal condition $ x(0)= \Gamma(x)$ can be more descriptive in physics with better effect than the classical
initial condition $x(0)=x_0,$ (see, e.g. \textup{\cite{Byszewski91,Byszewski95,Deng,Djebali}}).
 In the last case, i.e. $x(0)=x_0,$
 the problem (\ref{aplicationmulti1a}) has been studied by  Dhage \cite{Dhage 2005} and O'Regan \cite{Oregan} for the existence of solutions. Therefore it is of interest to discuss and to approximate
the solution of (\ref{aplicationmulti1a}) with a nonlocal initial condition  for various aspects of its solution under some suitable conditions.\\
 Similarly another class of nonlinear equations is used frequently to describe many phenomen a in various fields of applied sciences
 such as  physics, control theory, chemistry, biology, and so forth (see \cite{ASA2013}, \cite{Dhage 2006}, \cite{Jeribi2009} and \cite{JK2015}). This class is generated by the equations of the form:
  \begin{equation}\label{aplicationmultia}
   x(t)=  f(t,x(\sigma(t)))\cdot \left[q(t)+ \displaystyle\int_{0}^{\eta(t)}K(t,s,x(\tau(s)))ds\right], t\in J. \tag{P2}
\end{equation}
Both, (\ref{aplicationmulti1a}) and (\ref{aplicationmultia}), can be interpreted as fixed point problems in which the equation involved is a hybrid equation of the type
\begin{equation}\label{HE}x=Ax\cdot Bx. \tag{P3}\end{equation}
 A hybrid fixed
point theorem to (\ref{HE}) was proved by Dhage in \cite{Dhage 88}
 and since then, several extensions and generalizations of this fixed
point result have been proved. See \textup{\cite{Dhage05, Dhage06}} and the references therein. These results can be used to establish the existence and uniqueness of solutions. Although
the explicit calculation of the fixed point  is only possible in some
simple cases,  these  results are regarded as one of the most powerful tools  to approximate this fixed point by a computational method and to develop  numerical methods that allow  us  to approximate the solution of these equations.

 In Banach spaces, several works deals to develop numerical
techniques to approximate the solutions of some systems of integro-differential equations, by using different methods  such as the
Chebyshev polynomial method \cite{Sezer}, the parameterization method \cite{Dzhumabaev}, the wavelet methods \cite{Heydari},  the secant-like methods \cite{Argyros},
a collocation method in combination with operational matrices of Berstein polynomials \cite{Maleknejad}, the variational iteration
method \cite{Saberi}, etc. A combination method of a fixed point result and Schauder's basis in a
Banach space have been used in \textup{\cite{MI2010V, MI2013, MI2017, MI2020}} to solve numerically systems of integral and integro-differential equations.
 The advantages of this method over other numerical methods is its simplicity
  to implement it in a computer and the approximating functions are the sum of integrals of piecewise univariate and bivariate polynomials with coefficients easy to calculate.

Since the Banach algebras represents a practical framework for   equations such as   (\ref{aplicationmulti1a}) and (\ref{aplicationmultia}), and in general (\ref{HE}), the purposes of this paper are twofold. Firstly,    to present, under suitable conditions,  a method to approximate the fixed point of a hybrid equation of type (\ref{HE}), by means of  the product and composition of operators defined in a Banach algebra.  Secondly,  to develop and apply  the method presented  to  obtain an approximation of the solutions of  (\ref{aplicationmulti1a}) and (\ref{aplicationmultia}).

The structure of this paper is as follows:   in section \ref{sec:tools} we present some definitions and  auxiliary  results which will be needed in the sequel; in section \ref{sec:existence} we derive an approximation method for  the fixed point of the hybrid equation (\ref{HE}); in sections \ref{sec:dif} and \ref{sec:int}, we apply our results to prove the existence and the uniqueness of solution  of  (\ref{aplicationmulti1a}) and (\ref{aplicationmultia}),  we give an approximation method  for  these solutions and moreover, we establish some numerical examples to illustrate the applicability of our results.
\section{\textbf{Analytical tools}}\label{sec:tools}

 In this section, we  provide some  concepts and results that we  will  need  in  the  following  sections. The first analytical tool to be used comes from the theory of the fixed point.

\begin{defin}
 A mapping $A : X \longrightarrow X$ is said to be $\mathcal D$-Lipschitzian, if there exists a continuous nondecreasing
function $\phi: \mathbb{R}_+\longrightarrow  \mathbb{R}_+$ such that $$\|Ax-Ay\|\leq \phi(\|x-y\|)$$ for all $x,y \in X$, with $\phi(0)=0$. The mapping $\phi$ is called the $\mathcal D$-function associate to $A$. If $\phi(r)<r$ for $r > 0,$ the mapping $A$ is called a nonlinear contraction on $X$. $\hfill\diamondsuit$
\end{defin}
\begin{rem}\label{6}
 The class of $\mathcal D$-Lipschitzian mapping on $X$ contains the class of Lipschitzian mapping on $X$, indeed if $\phi(r)=\alpha\, r$, for some $\alpha>0$, then $A$ is called Lipschitzian mapping  with Lipschitz constant $\alpha$ or an $\alpha$-Lipschitzian mapping.  When $0\leq \alpha<1,$ we say that $A$ is a contraction. $\hfill\diamondsuit$
\end{rem}
The Banach fixed point theorem ensures that every contraction operator $A$ on a complete metric space $X$ has a unique fixed point $\tilde{x}\in X,$ and the sequence $\{A^nx, n\in \mathbb{N}\}$ converges to $\tilde{x}.$  One of the more useful generalizations of the Banach fixed point principal is the
following result due to Boyd and Wong in \cite{Boyd}.
\begin{thm}\label{Boyd-Wong}
 Let $(X, d)$ be a complete metric space, and let $A :
X \to X.$ Assume that  there exists a continuous function $\varphi : [0,\infty) \to [0,\infty)$ such that
$\varphi(r) < r$ if $r > 0,$ and
$$d( A (x), A (y)) \leq \varphi(d(x, y)), \ \ \forall x, y \in X.$$
Then $A$ has a unique fixed point $\tilde{x}\in X.$ Moreover, for any $x \in X,$ the sequence
$x_n = A^n(x)$ converges to $\tilde{x}.$
\end{thm}
\begin{rem}
The operator $A(x) = x - x^2$ mapping $X = [0, 1]$ into itself, and
possessing the unique fixed point $x = 0,$ does not satisfy the assumptions of the contraction principal (the smallest possible $k$ equals $1$), whereas it satisfies those of Theorem \ref{Boyd-Wong} with
$\varphi(r) = r - r^2.$
\end{rem}

On the other hand, Schauder bases will constitute the second essential tool.

A biorthogonal system in a Banach space $E$ is a system $\{(\tau_n, \xi_n), n\geq 1\}$ of $E\times E^*,$ where $E^*$ denotes the topological dual space of $E.$
Moreover, $\{(\tau_n, \xi_n), n\geq 1\}$ said to be a fundamental biorthogonal
system if $\overline{span}\{\tau_n\} = E.$
Now, a sequence $\{\tau_n, n\in \mathbb{N}\}\subset E$ defines a
Schauder basis of $E$ if, for every $x\in E,$ there is a unique sequence $(a_n)_n\subset \mathbb{R}$ such that $$x =\sum_{n\geq 1} a_n\tau_n.$$
This  produces the concept of the canonical sequence of finite dimensional projections $P_n : E \to E,$ defined by the formula
$$P_n\left(\sum_{k\geq 1} a_k\tau_k\right)= \sum_{k= 1}^n a_k\tau_k,$$ and the associated sequence of coordinate functionals $\tau_n^*\in E^*$  defined by the formula $$\tau^*_n\left(\sum_{k\geq 1} a_k\tau_k\right)= a_n.$$
Note that a Schauder basis is always a fundamental biorthogonal system, under the
interpretation of the coordinate functionals as biorthogonal functionals.
Moreover, in view of the Baire category theorem \cite{Brezis}, that for all $n\geq 1,$ $\tau_n^*$ and $P_n$ are continuous.
This yields, in particular, that
$$\lim_{n \rightarrow \infty}\|P_n(x)-x\|=0.$$
The above mentioned notions  play  an important  role to  approximating the solution of different integral and integro-differential equations (see \textup{\cite{MI2010V,MI2013,MI2017}}.)
\section{The  existence, uniqueness and approximation of a fixed point of the product of two operators in Banach algebras.}\label{sec:existence}
Based on  the  Boyd-Wong theorem, we establish the following fixed point result for the product of two nonlinear operators in Banach algebras.
\begin{thm}\label{eq}
Let $X$ be a nonempty closed convex subset of a Banach algebra $E.$ Let $A, B: X\to E$ be two operators  satisfying the following conditions:\\
\noindent $(i)$ $A$ and $B$ are $\mathcal{D}$-lipschitzian with $\mathcal D$-functions $\varphi$ and $\psi$ respectively,\\
\noindent $(ii)$ $A(X)$ and $B(X)$ are bounded,\\
\noindent $(iii)$ $Ax\cdot Bx\in X,$ for all $x\in X.$\\
\noindent Then, there is a unique point $\tilde{x}\in X$ such that $A\tilde{x} \cdot B\tilde{x}=\tilde{x}$ and the sequence $\{(A\cdot B)^nx\}, x \in X,$ converges to $\tilde{x}$ provided that $\|A(X)\|\psi(r)+\|B(X)\|\varphi(r)<r, r>0.$
$\hfill\diamondsuit$
\end{thm}
\noindent\textbf{Proof.}
Let $x, y\in X.$ we have
$$\begin{array}{rcl}\displaystyle\|Ax\cdot Bx- Ay\cdot By\|
 &\leq& \|Ax\cdot (Bx-By)\|+\|(Ax- Ay)\cdot By\|\\\\
 &\leq& \|Ax\|\,\|Bx-By\|+\|By\|\,\|Ax- Ay\|\\\\
  &\leq& \|A(X)\|\,\psi(\|x-y\|)+\|B(X)\|\,\varphi(\|x-y\|).
\end{array}
 $$
This implies that $A\cdot B$ defines a nonlinear contraction with $\mathcal{D}$-function
$$\phi(r)=\|A(X)\|\,\psi(r)+\|B(X)\|\,\varphi(r), \ r>0.$$
Applying the Boyd-Wong fixed point theorem \cite{Boyd}, we obtain the desired result. \hfill $\Box$\par\medskip

 \bigskip

Boyd-Wong's fixed-point theorem allows us to express the fixed point of $A\cdot B$ as the limit
of the sequence of functions $\{(A\cdot B)^n(x), n\in \mathbb{N}\},$ with
$x\in X.$ Obviously, if it were possible
to explicitly calculate, for each iteration, the expression $(A\cdot B)^n(x),$ then for each $n$ we
would have an approximation of the fixed point. But, as a practical matter, such an
explicit calculation is only possible in very particular cases. For this reason, we need to construct  another  approximation of the fixed point  which is   simple to calculate in practice.
Therefore, we need the following Lemmas. The proofs of these  Lemmas are similar to those of Lemma 1 and Lemma 2 in \cite{MI2020}.
\begin{lem}\label{lemma 1}
Let $X$ be a nonempty closed convex subset of a Banach algebra $E$ and let $A, B: X\to E$ and $\phi : \mathbb{R}^+ \to \mathbb{R}^+$ be a continuous nondecreasing function such that for all $n\geq 1,$
$$\left\|(A\cdot B)^nx-(A\cdot B)^ny\right\|\leq \phi^n (\|x-y\|).$$
Let
$x\in X$ and $T_0,T_1, \ldots, T_m: E\to E,$  with $T_0\equiv I.$  Then
$$\begin{array}{rcl}\left\|(A\cdot B)^mx-T_m\circ\ldots\circ T_1 x\right\|&\leq & \displaystyle\sum_{p=1}^{m-1}\phi^{m-p}\left(\left\|A\cdot B\circ T_{p-1}\circ\ldots\circ T_1 x-T_p\circ\ldots\circ T_1 x\right\|\right)\\\\
&&+\left\|A\cdot B\circ T_{m-1}\circ\ldots\circ T_1 x-T_m\circ\ldots\circ T_1 x\right\|.\end{array}$$ $\hfill\diamondsuit$
\end{lem}

\begin{lem}\label{lemma 2}
Let $X$ be a nonempty closed convex subset of a Banach algebra $E.$ Let $A, B: X\to E$ be two $\mathcal{D}$-Lipschitzian operators with $\mathcal{D}$-functions $\varphi$ and $\psi,$ respectively, and $A\cdot B$ maps $X$ into $X.$ Moreover, suppose
  that $$\phi(r):=\|A(X)\|\psi(r)+\|B(X)\|\varphi(r)<r, r>0.$$ Let $\tilde{x}$ be the unique fixed point of $A\cdot B,$ $x\in X,$ $\varepsilon>0,$ $n\in \mathbb{N}$ such that
$$\left\|(A\cdot B)^nx-T_n\circ\ldots \circ T_1 x\right\|\leq \frac{\varepsilon}{2},$$ then
$$\left\|\tilde{x}-T_n\circ\ldots \circ T_1 x\right\|\leq \varepsilon.$$ $\hfill\diamondsuit$
\end{lem}

Taking into account the above Lemmas, so as to approximate the solutions of the problems (\ref{aplicationmulti1a}) and (\ref{aplicationmultia}),  we will begin with an initial function
$x_0\in X$ and
we will construct a sequence of operators $\left\{S_n, n\in \mathbb{N}\right\}$ in order to obtain successive $T_n\circ \ldots \circ T_1(x_0)$
 approximations of the fixed point $\tilde{x}$ of the product
$A\cdot B$ following the scheme:

$$\begin{array}{ccc}
  x_0 &  &  \\
  \downarrow &  &  \\
  (A\cdot B)(x_0) & \approx & T_1(x_0)=A(x_0)\cdot S_1(x_0) \\
  \downarrow &  & \downarrow\\
  (A\cdot B)^2(x_0) & \approx & T_2\circ T_1x_0=(A\cdot S_2)\circ T_1(x_0) \\
  \vdots & \vdots & \vdots \\
  \\
  \vdots & \vdots & \vdots \\
   \downarrow&   & \downarrow\\
  (A\cdot B)^n(x_0) & \approx & T_n\circ\ldots\circ T_1(x_0)=(A\cdot S_n)\circ T_{n-1}\circ\ldots \circ T_1(x_0)\approx \tilde{x}
\end{array}$$

\section{\textbf{Nonlinear differential problems (\ref{aplicationmulti1a})}}\label{sec:dif}
In this section we focus  our  attention in the following  nonlinear differential equation with a nonlocal initial condition:
\begin{equation}\label{aplicationmulti1}\left\{
	\begin{array}{rl}
		\displaystyle \frac{d}{dt}\left(\frac{x(t)}{f(t,x(t))}\right)&=  g(t,x( t)), t\in J,\\\\
		x(0)&= \Gamma(x),
	\end{array} \tag{P1}
	\right.\end{equation}
where $J:=[0, \rho],$ $f: J\times \mathbb{R}\to \mathbb{R}\setminus\{0\},$ $g: J\times \mathbb{R}\to \mathbb{R}$ and $\Gamma : C(J) \to \mathbb{R}.$ Here $C(J)$ is the space  of all continuous functions from $J$ into $\mathbb{R}$ endowed with the norm $\|\cdot\|_\infty= \sup_{t\in J}|x(t)|.$ \\
This equation will be studied under the following assumptions:\\
\hspace*{20pt}{$(i)$} The partial mappings $t\mapsto f(t,x),$ $t\mapsto g(t,x)$ are continuous and the mapping $\Gamma$ is $L_\Gamma$-Lipschitzian.
\\
\hspace*{20pt}{$(ii)$} There exist $r>0$ and two nondecreasing, continuous functions $\varphi, \psi: \mathbb{R}_+ \longrightarrow \mathbb{R}_+ $ such that
$$\left|f(t,x)-f(t,{y})\right| \leq \alpha(t) \varphi(|x-y|),  t \in J, \hbox{ and }x, y \in \mathbb{R} \hbox{ with }|x|, |y|\leq r,$$
and
$$ \left|g(t,x)-g(t,{y})\right| \leq \gamma(t)\psi(|x-y|), t \in J \text{ and }x, y \in \mathbb{R}  \hbox{ with }|x|, |y|\leq r.$$
\hspace*{20pt}{$(iii)$} There is a constant $\delta>0$ such that $\sup_{x\in \mathbb{R}, |x| \leq r}|f(0,x)|^{-1}\leq \delta.$

\subsection{The existence and uniqueness of a solution to problem (\ref{aplicationmulti1}).}

In this subsection, we prove the existence and the uniqueness of a solution to the functional differential problem (\ref{aplicationmulti1}).
\begin{thm}\label{thm1}
Assume that the assumptions $(i)$-$(iii)$ hold. If
\begin{eqnarray*}\label{eq r011}
 \displaystyle\left\{
  \begin{array}{lll}
 \displaystyle M_F \delta  L_{\Gamma} t+\left(M_F\delta^2\alpha(0)\left(L_\Gamma r+\Gamma(0)\right)+M_G\|\alpha\|_\infty\right)\varphi(t)+M_F\|\gamma(\cdot)\|_{L^1} \psi(t)<t, \ t>0,\\\\
 \displaystyle  M_F M_G\leq r,
  \end{array}
\right.
\end{eqnarray*}
where $r$ is defined in the assumption $(ii),$ then the nonlinear differential problem  \eqref{aplicationmulti1}  has a unique solution in ${B}_r.$
\end{thm}
{\it Proof.} Let $$\Omega := \{x\in C(J); \|x\|\leq r\}.$$
 Here the constant $r$ is defined in $(ii).$
 Observe that $\Omega$ is a non-empty, closed, convex and bounded subset of $C(J),$ and the problem of the existence of a solution to  (\ref{aplicationmulti1}) can be formulated in the following fixed point problem $Fx\cdot Gx = x,$ where
$F, G $ are given for $x\in C(J)$ by
\begin{eqnarray}\label{ED1}\left\{
  \begin{array}{ll}
   (Fx)(t)&= \displaystyle f(t,x(t))\\\\
   (Gx)(t)&= \left[\displaystyle\frac{1}{f(0,x(0))}\Gamma(x)+ \displaystyle\int_{0}^{t}g(s,x(s))ds\right], t\in J.
 \end{array}
\right.\end{eqnarray}
 Let $x\in \Omega$ and $t, t'\in J.$  Since $f$ is $\mathcal{D}$-lipschitzian with respect to the second variable and is continuous with respect to the first variable, then by using the  inequality

 $$\begin{array}{rcl}\displaystyle|f(t,x(t))-  f(t',x(t'))|
 &\leq&  \displaystyle|f(t,x(t))-  f(t',x(t))|+|f(t',x(t))-  f(t', x(t'))|,
 \end{array}$$
we can show that $F$ maps $\Omega$ into $C(J).$ \\
Now, let us claim that $G$ maps $\Omega$ into $C(J).$ In fact,
  let $x\in \Omega$ and $t, t'\in J$ be arbitrary. Taking into account that $t\mapsto g(t,x)$ is a continuous mapping, it follows from assumption $(ii)$ that
 $$\begin{array}{rcl}\displaystyle| G(x)(t)-  G(x)(t')|
 &\leq&  \displaystyle\int_{t'}^{t}|g(s,x(s))-g(s,0)|ds+(t-t')\|g(\cdot,0)\|_\infty
\\\\
 &\leq& \displaystyle (t-t')\left(\|\gamma\|_{\infty}\psi(r)+\|g(\cdot,0)\|_\infty\right).
 \end{array}$$
This  proves  the claim.
Our strategy is to apply Theorem \ref{eq} to show the existence and the uniqueness of a fixed point for the product
$F\cdot G$ in $\Omega$ which in turn is a continuous solution for problem (\ref{aplicationmulti1}).

For this purpose, we will claim, first, that $F$ and $G$ are $\mathcal{D}$-lipschitzian mappings on $\Omega.$
 The claim regarding $F$ is clear in view of assumption  $(ii),$  that is
   $F$ is $\mathcal{D}$-lipschitzian with $\mathcal{D}$-function $\Phi$ such that
    $$\Phi(t)= \|\alpha\|_\infty\varphi(t), t\in J.$$
  We corroborate now the claim for $G.$ Let $x, y \in \Omega,$ and let $t\in J.$
 By using our assumptions, we obtain
 $$\begin{array}{rcl}\displaystyle\left|G(x)(t)-  G(y)(t)\right|
 &=& \left|\displaystyle\frac{1}{f(0,x(0))}\Gamma(x)-\frac{1}{f(0,y(0))}\Gamma(y)+
 \int_{0}^{t}g(s,x(s))-g(s,y(s))ds\right|\\\\
 &\leq&  \displaystyle \frac{L_{\Gamma}}{|f(0,x(0))|}\|x-y\|+\frac{\alpha(0)}{|f(0,x(0)) f(0,y(0))|}\left(L_\Gamma r+\Gamma(0)\right) \varphi(\|x-y\|)
 \\\\
&& +
 \displaystyle\int_{0}^{t}|\gamma(s)|\psi(|x(s)-y(s))|ds\\\\
 &\leq& \delta L_{\Gamma}\|x-y\|+\delta^2\alpha(0)\left(L_\Gamma r+\Gamma(0)\right)\varphi(\|x-y\|)+\|\gamma(\cdot)\|_{L^1} \psi(\|x-y\|).
 \end{array}$$
Taking the supremum over $t,$ we obtain that
  $G$ is $\mathcal{D}$-lipschitzian with $\mathcal{D}$-function $\Psi$ such that
   $$\Psi(t)=\delta  L_{\Gamma} t+\delta^2\alpha(0)\left(L_\Gamma r+\Gamma(0)\right)\varphi(t)+\|\gamma(\cdot)\|_{L^1} \psi(t), t\in J.$$
 On the other hand,  bearing  in mind assumption $(i),$ by using the above discussion we can see that $F(\Omega)$ and $G(\Omega)$
 are bounded  with bounds   $$M_F:=\|\alpha\|_\infty \varphi(r)+\|f(\cdot, 0)\|_\infty\hbox{ and } M_G:= \delta(L_\Gamma r+|\Gamma(0)|)+ \|\gamma\|_{\infty}\rho\psi(r)+\rho \|g(\cdot,0)\|_\infty.$$
 Taking into account the estimate $M_F M_G\leq r,$
   we obtain that $F\cdot G$ maps $\Omega$ into $\Omega.$ \\
   Now,  applying Theorem \ref{eq}, we infer that (\ref{aplicationmulti1})
has one and only one solution $\tilde{x}$ in $\Omega,$ and for each $x\in \Omega$ we have
$$\displaystyle\lim_{n\rightarrow \infty}(F\cdot G)^nx=\tilde{x}.$$
\hfill $\Box$\par\medskip

\noindent Notice that by induction argument we can show that
 \begin{eqnarray}\label{1a}\displaystyle\|(F\cdot G)^nx- (F\cdot G)^ny\|\leq \Theta^n(\|x-y\|),\end{eqnarray}
where $\Theta(t):=  M_F \Psi(t) + M_G \Phi(t), t\geq0.$

\subsection{Numerical method to approximate the solution of (\ref{aplicationmulti1}).}

In this subsection we find a numerical approximation of the solution to the nonlinear equation (\ref{aplicationmulti1}) using a   Schauder basis in $C(J).$

First, let us consider a Schauder basis $\{\tau_n\}_{n\geq1}$ in $C(J)$ and the sequence of
associated projections $\{\xi_n\}_{n\geq 1}.$
Let  $$\left\{
                           \begin{array}{ll}
                              T_p: C(J)\longrightarrow C(J)\\\\
                             x\mapsto \displaystyle T_p(x)(t)
                              =F(x)(t) \left( \displaystyle\frac{1}{f(0,x(0))}
                             \Gamma(x)+ \displaystyle\int_{0}^{t}\xi_{n_p}(U_0(x))(s)ds\right),
                           \end{array}
                         \right.
$$
where $F: C(J)\longrightarrow C(J)$ such that $$F(x)(t)=f(t,x(t))$$
 and
 $U_0: C(J)\longrightarrow C(J)$ such that $$U_0(x)(s)=g(s,x(s)).$$
 \\
 \begin{rem}\label{rem1}
$(i)$ For all fixed $p\geq 1,$ the mapping $T_p$ maps $\Omega$ into $\Omega.$\\
 In fact, let $x\in \Omega,$ we have
 $$\begin{array}{rcl} \left|T_p(x)(t)\right|&=&
\left| F(x)(t)  \left( \displaystyle\frac{1}{f(0,x(0))}\Gamma(x)+ \displaystyle\int_{0}^{t}\xi_{n_p}(U_0(x))(s)ds\right|\right)\\\\
&\leq&
\left|f(t,x(t))\right|  \left( \displaystyle\delta |\Gamma(x)|+ \int_{0}^{t}\left|\xi_{n_p}(U_0(x))(s)\right|ds\right).
\end{array}$$
Proceeding essentially as in the above subsection and using the fact that $\xi_{n_p}$ is a bounded linear operator on $C(J),$ we get
$$\begin{array}{rcl} \left|T_p(x)(t)\right| &\leq&
M_F  \left[\displaystyle\delta |\Gamma(x)|+ \rho
\left\|\xi_{n_p}\left(U_0(x)\right)\right\|\right]\\\\
 &\leq&
M_F \displaystyle \left[\displaystyle\delta(L_\Gamma r+|\Gamma(0)|)+\rho \sup_{s\in J}|g(s,x(s))|\right]\\\\
 &\leq&
M_F M_G.
\end{array}$$
 In view of assumption $(iii),$ we infer that $T_p$ maps $\Omega$ into $\Omega.$
 \\

\noindent $(ii)$ Item $(i)$ means, in particular, that for all fixed $p\geq 1,$ the operator  $T_p\circ\ldots\circ T_1$ maps $\Omega$ into $\Omega.$$\hfill\diamondsuit$
 \end{rem}

 Our objective is to justify that we can choose $n_1, \ldots, n_m,$ so that the operators
$T_{1}, \ldots ,T_{m}$ can be used to obtain an approximation to the unique solution of equation (\ref{aplicationmulti1}).
\begin{thm} Let $\tilde{x}$ be the unique solution to the nonlinear problem (\ref{aplicationmulti1}). Let $x\in \Omega$ and $\varepsilon>0,$ then there exists $n\in \mathbb{N}$ such that
$$\left\|\tilde{x}-T_n\circ\ldots \circ T_1 x\right\|\leq \varepsilon.$$ $\hfill\diamondsuit$
\end{thm}
\noindent{\it Proof.} Let $x\in \Omega$ and $\varepsilon>0.$ For $p \in \{1, \ldots ,m\},$ we define $U_p : C(J) \to C(J)$ by
$$U_p(x)(s) := g(s,T_p \circ \ldots \circ T_1(x)( s)),
 \ s\in J, x\in C(J)$$
 and
  $F_p : C(J) \to C(J)$ by
$$F_p(x)(s) := F\left(s,T_p \circ \ldots \circ T_1(x)( s)\right),
 \ s\in J, x\in C(J).$$
 According to  inequality (\ref{1a}), in view of Lemma \ref{lemma 1}, we get $$ \left\|(F\cdot G)^mx-T_m\circ\ldots \circ T_1 x\right\| \leq $$  $$\displaystyle \sum_{p=1}^{m-1}\Theta^{m-p}\left(\left\|(F\cdot G)\circ T_{p-1}\circ\ldots\circ T_1 x-T_p\circ\ldots\circ T_1 x\right\|\right)
  +\left\|(F\cdot G)\circ T_{m-1}\circ\ldots\circ T_1 x-T_m\circ\ldots\circ T_1 x\right\|. $$
 Taking into account (\ref{rem1})-$(i)$ and using similar arguments as in the above section,  we infer that $\left\|F_{p-1}(x)\right\|$ is bounded, and consequently we get
 $$\begin{array}{rcl}&&\displaystyle\left|(F\cdot G)\circ T_{p-1}\circ\ldots\circ T_1(x)(t)- T_p\circ T_{p-1}\circ\ldots\circ T_1(x)(t)\right|\\\\
 &\leq& \left|F_{p-1}(x)(t)\left(\displaystyle\int_{0}^t g\left( s,T_{p-1}\circ\ldots\circ T_1(x)(s)\right)\, ds-\displaystyle\int_{0}^t \xi_{n_p}(U_{p-1}(x))(s)\, ds\right)\right|\\\\
 &\leq&
 \left|F_{p-1}(x)(t)\right|\, \displaystyle\int_{0}^t \left|\left(\xi_{n_p}(U_{p-1}(x))-U_{p-1}(x)\right)(s)\right|\, ds\\\\
 &\leq&
 \rho\left\|F_{p-1}(x)\right\| \, \left\|\xi_{n_p}(U_{p-1})(x)-U_{p-1}(x)\right\|
 .\end{array}$$
 Then, we obtain
 $$ \left\|(F\cdot G)^mx-T_m\circ\ldots \circ T_1 x\right\| \leq   \displaystyle\sum_{p=1}^{m-1}\Theta^{m-p}\left(
\rho M_F \, \left\|\xi_{n_p}(U_{p-1}(x))-U_{p-1}(x)\right\|\right)
 +
\rho M_F\, \left\|\xi_{n_m}(U_{m-1}(x))-U_{m-1}(x)\right\|. $$
In view of the convergence property of the Projection operators associated to the Schauder  basis and the continuity of $\Theta,$  we can find $n_1, \ldots , n_m \geq 1$ and
therefore $T_1, \ldots,T_m,$ such that
 $$ \|(F\cdot G)^mx-T_m\circ\ldots \circ T_1 x\| \leq$$   $$\displaystyle\sum_{p=1}^{m-1}\Theta^{m-p}\Big(
 \rho M_F \left\|\xi_{n_p}(U_{p-1}(x))-U_{p-1}(x)\right\|\Big)
 +\rho  M_F \left\|\xi_{n_m}(U_{m-1}(x))-U_{m-1}(x)\right\|
 \leq   \displaystyle \frac{\varepsilon}{2}. $$
Now apply Lemma \ref{lemma 2}, in order to get $$\left\|\tilde{x}-T_m\circ\ldots \circ T_1(x)\right\|<\varepsilon.$$
 \hfill $\Box$\par\medskip
\subsection{Numerical experiments.}

This section is devoted to  providing  some examples and their numerical results to illustrate the theorems of the above sections. We will consider $J=[0,1]$ and the classical Faber-Schauder system  in $C(J)$ where the nodes are  the naturally ordered dyadic numbers (see \textup{\cite{Gelbaum, Semadeni}} for details).

\begin{exam}
Consider the nonlinear differential equation with a nonlocal initial condition
\begin{equation}\label{Ex1}\left\{
    \begin{array}{ll}
&\displaystyle \frac{d}{dt}\left(\frac{x(t)}{f(t,x(t))}\right)=  \displaystyle  a e^{-x(t)}, \ \ t\in J,\\\\
  & \displaystyle x(0)= b( \sup_{t\in J}|x(t)|+3/4)
\end{array}\right.\end{equation}
\noindent where $f(t,x)=\displaystyle\frac{b}{1+ae^{-b}t},$ $g(t,x)=\displaystyle ae^{-x},$ $\Gamma(u) =b\left( \sup_{t\in J}|u(t)|+3/4\right),$ with  $a<1/\log(2).$
 \\
Let $R$ be small enough such that $a(\log(2)+R)<1.$ Let $x, y \in [-R,R],$ by an elementary calculus we can show that
$$\begin{array}{rcl}\left|f(t,x)-f(t,y)\right|&\leq & \displaystyle \alpha(t) \varphi(|x-y|) \end{array}$$
and $$\begin{array}{rcl}\left|g(t,x)-g(t,y)\right|&\leq & \displaystyle \gamma(t) \psi(|x-y|) \end{array}$$
where  $\displaystyle \alpha(t)=\varphi(t)=0,$ $ \gamma(t)=a e^R(1-e^{-t}),$ and $\psi(t)=t.$ \\
On the other hand, we have that
%
$\Gamma$ is Lipschizian with a Lipschiz constant $L_\Gamma=b,$ and
 $$\displaystyle\sup_{x,|x|\leq R}[f(0,x)]^{-1} \leq \delta=\frac{1}{b}.$$   Applying  Theorem \ref{thm1}, we obtain that  (\ref{Ex1}) has a unique solution in
 $\Omega=\left\{x\in C([0,1]); \|x\|\leq 3/4\right\},$ when $a$ is  small enough. In fact the solution is $x(t)=b.$ We apply the numerical method for $a=0.1$ and the initial $x_0(t)=\frac{1}{4}\left(\sqrt{bt}+ 1\right).$
\begin{center}
	\small{\textsc{Table 1. Numerical results for (\ref{Ex1}) with initial $x_0(t)=\frac{1}{4}\left(\sqrt{bt}+ 1\right)$.}}

	\begin{tabular}{|c|c|c|c|c|c| }
		\hline \hline
		& 	& \multicolumn{2}{c|}{$n_1=\dots=n_m=9 $ } & \multicolumn{2}{c|}{$n_1=\dots=n_m=33$ }  \\
		\hline
		$t$  &$x^*(t)$	& $m=2$ & $m= 4$   & $m= 2$ & $m= 4$    \\
		\hline \hline
		$0.1$& $0.25$	&0.25964068562641207 & 0.2526360625738145 &0.25779577744548676 & 0.25062384017038686 \\
		\hline
		$0.2$&  $0.25$ & 0.2581608132685836 & 0.2512245431325148 & 0.2576778369861067& 0.2506151528771704 \\
		\hline
		$0.3$ & $0.25$	& 0.25785705013803817
& 0.25102089532293176 &  0.2575623161293616& 0.25060665510642743\\
		\hline
		$0.4$& $0.25$	&   0.25774159101031774&
0.25100874582984495 &0.25744919245075903 & 0.2505983412941664\\
		\hline
		$0.5$& $0.25$	& 0.2576285346430475 &0.2509968386936278 &  0.25733843642963494&   0.2505902060799007\\
		\hline
		$0.6$ & $0.25$	&   0.25751784650251447 &  0.2509851672563384&   0.2572300145734221&0.2505822442972077\\
		\hline
		$0.7$ & $0.25$	&0.2574094899054699
&0.25097372508850474&0.2571238909612659&0.2505744509661791\\
		\hline
		$0.8$& $0.25$	& 0.25730342712624404 &0.2509625059364119&0.25702002815700664&  0.25056682128612107\\
		\hline
		$ 0.9$ &$ 0.25$	&  0.25719961932300417& 0.25095150376429876& 0.2569183878122034& 0.25055935062726176\\
			\hline
		$1 $   &  $0.25$ &0.2570980270442474&0.2509407127451644 &0.25681893107062653&0.2505520345235613\\
		\hline \hline
	\end{tabular}
\end{center}

\begin{center}
	\begin{tabular}{|c|c|c|c|c| }
		\hline \hline
		 	& \multicolumn{2}{c|}{$n_1=\dots=n_m=9 $ } & \multicolumn{2}{c|}{$n_1=\dots=n_m=33$ }  \\
		\hline
		  	& $m=2$ & $m= 4$   & $m= 2$ & $m= 4$    \\
		\hline \hline
		 $\|x^*-\tilde{x}\|_\infty$	&  $9.90603\times 10^{-3}$&  $2.86369\times 10^{-3}$  &$8.33966\times 10^{-3}$ &  $1.0862\times 10^{-3}$  \\
		\hline \hline
	\end{tabular}
\end{center}
\end{exam}

\begin{exam}
Consider the nonlinear differential equation with a nonlocal initial condition
\begin{equation}\label{Ex2}\left\{
    \begin{array}{ll}
&\displaystyle \frac{d}{dt}\left(\frac{x(t)}{f(t,x(t))}\right)=  \displaystyle  a(x(t))^2, \ \ t\in J,\\\\
 &\displaystyle x(0)= 1/(4b)\sup_{t\in J}|x(t)|^2,
    \end{array}
  \right.
\end{equation}
where $f(t,x)=\displaystyle\frac{b(t+1)}{1+\frac{a b^2}{3}(x^3/b^3-1)}$ and $g(t,x)=a x^2,$ with $ab^2<3.$
 \\
Let $R>0$ such that $2b\leq R$ and $\frac{a}{3b}(b^3+R^3)<1.$ Let $x, y \in[-R,R].$ By an elementary calculus we can show that
$$\begin{array}{rcl}\left|f(t,x)-f(t,y)\right|&\leq & \displaystyle \alpha(t) \varphi(|x-y|) \end{array}$$
and $$\begin{array}{rcl}\left|g(t,x)-g(t,y)\right|&\leq & \displaystyle \gamma(t) \psi(|x-y|) \end{array}$$
where  $ \displaystyle\alpha(t)=\frac{a(t+1) R^2}{\left(1-\frac{a}{3b}(R^3+b^3)\right)^2},$ $ \gamma(t)=2a R,$ and $\varphi(t)=\psi(t)=t.$ \\
On the other hand, we have that
$$\begin{array}{rcl}\displaystyle|\Gamma(u)-\Gamma(v)|&\leq & \displaystyle  \frac{R}{2b}\|u-v\|.
 \end{array}$$
Consequently, $\Gamma$ is Lipschizian with aLipschiz constant $L_\Gamma=\frac{R}{2b}.$
 It is easy to prove that
$$\sup_{x\in \mathbb{R},|x|\leq R}[f(0,x)]^{-1} \leq \delta=aR^3/(3b^2)+1/b.$$
 Now,  applying Theorem \ref{thm1}, in order to obtain that   (\ref{Ex2}), with $a$ is
  small enough, has a unique solution in
$\Omega=\left\{x\in C([0,1]); \|x\|\leq 1/2\right\}.$ We can check that the solution is $x(t)=b(t+1).$

\noindent The following table  shows  the numerical results of the proposed method for $a=0.05$, $b=1/4$ and $x_0(t)=\frac{1}{2} t.$
\begin{center}
	\small{\textsc{Table 2. Numerical results for    (\ref{Ex2})   with initial $x_0(t)=\frac{1}{2} t$.}}

\begin{tabular}{|c|c|c|c|c|c| }
		\hline \hline
		& 	& \multicolumn{2}{c|}{$n_1=\dots=n_m=9 $ } & \multicolumn{2}{c|}{$n_1=\dots=n_m=33$ }  \\
		\hline
		$t$  &$x^*(t)$	& $m=2$ & $m= 4$   & $m= 2$ & $m= 4$    \\
		\hline \hline
		$0.1$& $0.275$	&0.27417118067351837& 0.27151545133640886 & 0.2740819659311709 & 0.27145329704728827 \\
		\hline
		$0.2$&  $0.3$	      & 0.2990105059004299 & 0.29611673530305527 & 0.2989981055587191&   0.29613324654650613 \\
		\hline
		$0.3$ & $0.325$	& 0.32391591487977106&  0.3207837845940706 & 0.3239149335622202  & 0.3208140511167786 \\
		\hline
		$0.4$& $0.35$	&
0.3488342313971621 &
0.34546352791535867 &0.3488329357145739&0.34549585475483185 \\
		\hline
		$0.5$& $0.375$	& 0.3737541821287371&  0.3701445199310059&  0.3737524860894689  &  0.3701788114857308 \\
		\hline
		$0.6$ & $0.40$	& 0.39867580493805804 &  0.3948268789541488 & 0.39867366683899425&  0.39486308640853285 \\
		\hline
		$0.7$ & $0.425$	&
0.42359867388137695  &
0.4195107187398104 &0.4235961144223516 &0.41954885401447617 \\
		\hline
		$0.8$& $0.45$	&  0.4485219829977168& 0.4441962543294659 & 0.448518998226605&  0.44423629583080837  \\
		\hline
		$ 0.9$ &$ 0.475$	& 0.47344474190449987 & 0.4688837174935067 &0.47344130014978747 &  0.468925600958778\\
			\hline
		$1 $   &  $ 0.5$  &
0.498366567564148 &0.49357335586512446  &0.49836264112050677 & 0.4936169655580174 \\
		\hline \hline
	\end{tabular}
\end{center}

\begin{center}
	\begin{tabular}{|c|c|c|c|c| }
		\hline \hline
		 	& \multicolumn{2}{c|}{$n_1=\dots=n_m=9 $ } & \multicolumn{2}{c|}{$n_1=\dots=n_m=33$ }  \\
		\hline
		  	& $m=2$ & $m= 4$   & $m= 2$ & $m= 4$    \\
		\hline \hline
		 $\|x^*-\tilde{x}\|_\infty$	& $1.63343\times 10^{-3}$&$6.42664\times 10^{-3}$ &  $1.63736\times 10^{-3}$ & $6.38303\times 10^{-3}$ \\
		\hline \hline
	\end{tabular}
\end{center}
\end{exam}

\section{\textbf{Nonlinear integral equations of type (\ref{aplicationmultia})}.}\label{sec:int}
This section  deals with the following nonlinear integral equation:
\begin{equation}\label{aplicationmulti}
   x(t)=  f(t,x(\sigma(t)))\cdot \left[q(t)+ \displaystyle\int_{0}^{\eta(t)}K(t,s,x(\tau(s)))ds\right], t\in J.\tag{P2}
\end{equation}
where $\sigma, \tau, \eta: J\to J,$ $f: J\times \mathbb{R}\to \mathbb{R}, q\in C(J)$ and $K: J\times J\times \mathbb{R}\to \mathbb{R}.$\\
More precisely, we prove the existence and the uniqueness of a solution to equation (\ref{aplicationmulti}), and then we  provide  an approximation method of this solution.
In our consideration, we need the following hypotheses:\\
\hspace*{20pt}{$(i)$} The partial mappings $t\mapsto f(t,x)$ and $(t,s)\mapsto K(t,s,x)$ are continuous.
\\
\hspace*{20pt}{$(ii)$} There exist $r>0$ and two nondecreasing, continuous functions $\varphi, \psi: \mathbb{R}_+ \longrightarrow \mathbb{R}_+ $ such that
$$ \left|f(t,x)-f(t,{y})\right| \leq \alpha(t) \varphi(|x-y|),  t \in J, \hbox{ and }x, y \in \mathbb{R} \hbox{ with }|x|, |y|\leq r,$$
and
$$ \left|K(t,s,x)-K(t,s,{y})\right| \leq \gamma(t,s)\psi(|x-y|), t, s \in J \text{ and }x, y \in \mathbb{R}  \hbox{ with }|x|, |y|\leq r.$$

\subsection{The existence and uniqueness of a solution to Eq. (\ref{aplicationmultia})}\label{5.1}
To allow the abstract formulation of equation (\ref{aplicationmulti}), we define
the following operators by
\begin{eqnarray}\label{EI}\left\{
  \begin{array}{ll}
   (Fx)(t)&= \displaystyle f(t,x(\sigma(t)))\\\\
   (Gx)(t)&= \left[\displaystyle q(t)+ \displaystyle\int_{0}^{\eta(t)}K(t,s,x(\tau(s)))ds\right], t\in J.
 \end{array}
\right.\end{eqnarray}
First, we will establish the following result which  shows  the existence and uniqueness  of a solution.
\begin{thm}\label{thm2}
Assume that  the assumptions $(i)$-$(ii)$ hold. If
\begin{eqnarray}\label{eq r011}
 \displaystyle\left\{
  \begin{array}{ll}
  \displaystyle M_F\rho \|\gamma\|_\infty\psi(t) + M_G \|\alpha\|_\infty \varphi(t)<t, \ t>0\\\\
 \displaystyle  M_F M_G\leq r,
  \end{array}
\right.
\end{eqnarray}
where $$M_F=\|\alpha\|_\infty \varphi(r)+\|f(\cdot,\theta)\|_\infty\hbox{ and } M_G=\displaystyle\|q(\cdot)\|_\infty+\rho\left(\|k(\cdot,\cdot,0)\|_\infty+
\|\gamma\|_\infty \psi(r)\right),$$ then the nonlinear integral equation \eqref{aplicationmulti}
 has a unique solution in ${B}_r.$
\end{thm}
{\it Proof.}
Let $\Omega := \{x\in C(J); \|x\|_\infty\leq r\}.$ By using similar arguments to those in  the above section, we can show that $F$  and $G$ define $\mathcal{D}$-lipschitzian mappings from $\Omega$ into $C(J),$ with $\mathcal{D}$-functions $\|\alpha\|_\infty\varphi$ and $\rho \|\gamma\|_\infty \psi,$ respectively.
 Also it is easy to see that $F(\Omega)$ and $G(\Omega)$
 are bounded with bounds, respectively,
 $$M_F=\|\alpha\|_\infty \varphi(r)+\|f(\cdot, 0)\|_\infty\hbox{ and }
  M_G=\|q\|_\infty
  + \rho\left(\|\gamma\|_\infty\psi(r)+\|K(\cdot,\cdot, 0)\|_\infty\right).$$
 Taking into account our assumptions, we deduce that $F\cdot G$ maps $\Omega$ into $\Omega.$ Now, an application of
Theorem \ref{eq} yields that (\ref{aplicationmulti})
has one and only one solution $\tilde{x}$ in $\Omega,$ and for each $x\in \Omega$ we have $\displaystyle\lim_{n\rightarrow \infty}(F\cdot G)^nx=\tilde{x}.$ \hfill $\Box$\par\medskip
\noindent Moreover, by induction argument we can obtain
 \begin{eqnarray}\label{1}\displaystyle\|(F\cdot G)^nx- (F\cdot G)^ny\|\leq \Theta^n(\|x-y\|),\end{eqnarray}
where $\Theta(t):=\rho \|\gamma\|_\infty M_F \psi(t) +\|\alpha\|_\infty M_G \varphi(t), t\geq0.$
\subsection{A Numerical method to approximate the solution to  (\ref{aplicationmulti}).}
In this subsection we  provide  a numerical approximation of the solution to the nonlinear equation (\ref{aplicationmulti}).
Now let us consider a Schauder basis $\{\tau_n\}_{n\geq1}$ in $C(J\times J)$ and the sequence of
associated projections $\{\xi_n\}_{n\geq 1}.$
Let  $$\left\{
                           \begin{array}{ll}
                              T_p: C(J)\longrightarrow C(J)\\\\
                             x\mapsto \displaystyle T_p(x)(t)
                             =F(x)(t) \left({q}(t)+ \displaystyle\int_{0}^{\eta(t)}\xi_{n_p}(U_0(x))(t,s)ds\right),
                           \end{array}
                         \right.
$$
where
  $F: C(J)\longrightarrow C(J)$ such that $$F(x)(t)=f(t,x(\sigma(t)))$$
 and
 $U_0: C(J)\longrightarrow C(J\times J)$ such that $$U_0(x)(t, s)=K(t,s,x(\tau(s))).$$

 \begin{rem}\label{rem2}
$(i)$ For all fixed $p\geq 1,$ the mapping $T_p$ maps $\Omega$ into $\Omega.$\\
 In fact, let $x\in \Omega,$ we have
 $$  \left|T_p(x)(t)\right| =
\left| F(x)(t) \left(\displaystyle {q}(t)+ \displaystyle\int_{0}^{\eta(t)}\xi_{n_p}(U_0(x))(t,s)ds\right)\right|\\ \leq
\left|f(t,x(\sigma(t)))\right| \,  \left(\displaystyle \left|{q}(t)\right|+ \displaystyle\int_{0}^{\eta(t)}\left|\xi_{n_p}(U_0(x))(t,s)\right|ds\right). $$
Proceeding essentially  as in the above section and using the fact that $\xi_{n_p}$ is a bounded linear operator on $C(J\times J),$ we get
$$ \left|T_p(x)(t)\right|  \leq
M_F  \left(\left|q(t)\right|+ \displaystyle \rho
\left\|\xi_{n_p}\left(U_0(x)\right)\right\|\right)
  \leq
M_F  \left(\left\|q\right\|_\infty+ \displaystyle \rho \sup_{t,s\in J}|k(t,s,x(\tau(s)))|\right)\\\\
  \leq
M_F  M_G.
 $$
 In view of our assumptions, we infer that $T_p$ maps $\Omega$ into $\Omega.$
 \\

\noindent $(ii)$ Item $(i)$ means, in particular, that for all fixed $p\geq 1,$ the operator  $T_p\circ\ldots\circ T_1$ maps $\Omega$ into $\Omega.$
 \end{rem}

 Again, our objective is to justify that we can choose $n_1, \ldots, n_m,$ so that the operators
$T_{1}, \ldots ,T_{m}$ can be used to obtain an approximation of the unique solution to equation (\ref{aplicationmulti}).
\begin{thm} Let $\tilde{x}$ be the unique solution to the nonlinear equation (\ref{aplicationmulti}). Let $x\in \Omega$ and $\varepsilon>0,$ then there exists $n\in \mathbb{N}$ such that
$$\left\|\tilde{x}-T_n\circ\ldots \circ T_1 x\right\|\leq \varepsilon.$$ $\hfill\diamondsuit$
\end{thm}
\noindent{\it Proof.} Let $x\in \Omega$ and $\varepsilon>0.$ For $p \in \{1, \ldots ,m\},$ we define $U_p : C(J) \to C(J\times J)$ by
$$U_p(x)(t,s) := K(t,s,T_p \circ \ldots \circ T_1(x)(s)),
 \ t, s\in J, x\in C(J)$$
 and
  $F_p : C(J) \to C(J)$ by
$$F_p(x)(s) := f\left(s,T_p \circ \ldots \circ T_1(x)(s)\right),
 \ s\in J, x\in C(J).$$
 According to Lemma \ref{lemma 1}, we get $$ \left\|(F\cdot G)^mx-T_m\circ\ldots \circ T_1 x\right\| \leq $$ $$ \displaystyle \sum_{p=1}^{m-1}\Theta^{m-p}\left(\left\|(F\cdot G)\circ T_{p-1}\circ\ldots\circ T_1 x-T_p\circ\ldots\circ T_1 x\right\|\right)
  +\left\|(F\cdot G)\circ T_{m-1}\circ\ldots\circ T_1 x-T_m\circ\ldots\circ T_1 x\right\|. $$
 Taking into account Remark \ref{rem2},  we infer that $\left\|F_{p-1}(x)\right\|$ is bounded.
 Proceeding essentially, as in the above section, we get
 $$\begin{array}{rcl}\displaystyle\left|(F\cdot G)\circ T_{p-1}\circ\ldots\circ T_1(x)(t)- T_p\circ T_{p-1}\circ\ldots\circ T_1(x)(t)\right|\leq
   \rho\left\|F_{p-1}(x)\right\| \, \left\|\xi_{n_p}(U_{p-1})(x)-U_{p-1}(x)\right\|,\end{array}$$
 which implies that
 $$ \left\|(F\cdot G)^mx-T_m\circ\ldots \circ T_1 x\right\| \leq  $$ $$\displaystyle\sum_{p=1}^{m-1}\Theta^{m-p}\left(
\rho M_F \, \left\|\xi_{n_p}(U_{p-1})(x)-U_{p-1}(x)\right\|\right)
 +
\rho M_F \, \left\|\xi_{n_m}(U_{m-1})(x)-U_{m-1}(x)\right\|. $$
In view of the convergence property of the Projection operators associated to the Schauder basis, we can find $n_1, \ldots , n_m \geq 1$ and
therefore $T_1, \ldots,T_m,$ such that
 $$ \|(F\cdot G)^mx-T_m\circ\ldots \circ T_1 x\| \leq $$ $$\displaystyle\sum_{p=1}^{m-1}\Theta^{m-p}\Big(
 \rho M_F \left\|\xi_{n_p}(U_{p-1})(x)-U_{p-1}(x)\right\|\Big)
 +
\rho  M_F \left\|\xi_{n_m}(U_{m-1})(x)-U_{m-1}(x)\right\|\\\\
 \leq  \displaystyle \frac{\varepsilon}{2}. $$
Now apply Lemma \ref{lemma 2}, in order to get $$\|\tilde{x}-T_m\circ\ldots \circ T_1(x)\|<\varepsilon.$$
\hfill $\Box$\par\medskip

\subsection{\textbf{Numerical experiments.}}

This section is devoted to give some examples and their numerical results to illustrate  the previous results   using the usual Schauder basis in $C([0,1]^2)$  with the well know square ordering (see for example \cite{Gelbaum}).
\begin{exam}
Consider the nonlinear differential equation
\begin{equation}\label{Ex3}
   \displaystyle x(t)=  \displaystyle a(t+1)\left[\frac{b}{a}-\frac{b^2}{3}\left((t+1)^3-1\right)+\int_0^t(x(s))^2 ds\right], \ \ \ t\in J.
\end{equation}
Proceeding essentially as in subsection \ref{5.1}, equation (\ref{Ex3}) can be written as  a fixed point problem $$x=F(x)\cdot G(x),$$ where $F$ and $G$ are defined in (\ref{EI}), with $f(t,x)=a (t+1),$ $q(t)=b/a-\frac{b^2}{3}\left((t+1)^3-1\right)$ and $k(t,s,x)=x^2.$\\
Let $x,y\in [-R,R],$ we have that
  $$\left|k(t,s,x)-k(t,s,y)\right|\leq\gamma(t,s)\psi(|x-y|)$$
where $\displaystyle \gamma(t,s)=2R,$ and $\psi(t)=t.$\\
 An application of Theorem \ref{thm2}, yields that   (\ref{Ex2}) has a unique solution in $\Omega=\left\{x\in C([0,1]); \|x\|\leq3\right\}$, in fact the solution is $x(t)=b(t+1).$
\\
Using the proposed method with $a=0.1,$ $b=0.1$ and $x_0(t)=t^2,$ we obtain the following table:
\begin{center}
\small{\textsc{Table 3. Numerical results for  the   (\ref{Ex3})}}  with initial $x_0(t)=t^2$.
\begin{tabular}{|c|c|c|c|c|c|}
	\hline \hline
	& 	& \multicolumn{2}{c|}{$n_1=\dots=n_m=9 $ } & \multicolumn{2}{c|}{$n_1=\dots=n_m=33  $ }  \\
	\hline
	$t$  &$x^*(t)$	& $m=2$ & $m=4$ & $m=2$ & $m=4$    \\
	\hline \hline
	$0.1$& $0.11$	& 0.10994463333333335 & 0.10994463333333335 & 0.10995952813954675 &  0.10995955765685321 \\
	\hline
	$0.2$&  $0.12$	&  0.11981787538985864&  0.11981791805770493&  0.11994695097301926 &   0.11994727822516114\\
	\hline
	$0.3$ &$0.13$	& 0.12975090261315447 & 0.12975116990203317& 0.12993149120265635& 0.1299327014013851 \\
	\hline
	$0.4$& $0.14$	&
0.1396858063225927&
0.13968664031615474 &
0.13991267664284926 &
0.13991561146443787\\
	\hline
	$0.5$& $0.15$	&  0.14960946152701812 &0.1496116012197044& 0.14989041600508018 & 0.14989578496520412 \\
	\hline
	$0.6$ & $0.16$	&  0.15952079194994145&  0.1595251486759711&  0.15986592032169128 & 0.15987299132148372 \\
	\hline
	$0.7$ & $0.17$	&
0.16941936215524034  &
0.1694262809122463&
0.1698435692932222&
0.1698469898893412 \\
	\hline
	$0.8$& $0.18$	& 0.17930673541360537 &0.1793140741901599& 0.17983433359550066 &  0.17981752625254807\\
	\hline
	$ 0.9$ &$ 0.19$	&   0.18918899833227526  &   0.18918756887790722
& 0.18986179093331174 &   0.1897843325246908 \\
	\hline
	$ 1$ &$ 0.2$	&
0.19908194969111695 &
0.1990457618518603&
0.1999725602822185 &   0.1997471266515799 \\
	
	\hline \hline
\end{tabular}
\end{center}

\begin{center}
	\begin{tabular}{|c|c|c|c|c| }
		\hline \hline
		 	& \multicolumn{2}{c|}{$n_1=\dots=n_m=9 $ } & \multicolumn{2}{c|}{$n_1=\dots=n_m=33$ }  \\
		\hline
		  	& $m=2$ & $m= 4$   & $m= 2$ & $m= 4$    \\
		\hline \hline
		 $\|x^*-\tilde{x}\|_\infty$	& $9.1805\times 10^{-4}$ & $9.544238\times 10^{-4}$  & $1.65588\times 10^{-4}$ &  $2.52873\times 10^{-4}$  \\
		\hline \hline
	\end{tabular}
\end{center}
\end{exam}

\begin{exam}
Consider the nonlinear differential equation
\begin{equation}\label{Ex4}
\displaystyle x(t)=  \displaystyle \left(a e^{-x(t)}+b\right)\left[\frac{t}{ae^{-t}+b}
+\frac{1}{1-c}\log(\cos(1-c)t)
+\int_0^t \tan(1-c)x(s) ds\right].
\end{equation}
 Similarly  to that above,   (\ref{Ex4}) can be written as a fixed point problem $$x=F(x)\cdot G(x),$$ with the same notations in (\ref{EI}).\\
Let $R>0$ and let $x, y\in[-R,R].$ By an elementary calculus we can show that
$$|f(t,x)-f(t,y)|\leq\alpha(t)\varphi(|x-y|)$$
and $$|k(t,s,x)-k(t,s,y)|\leq\gamma(t)\psi(|x-y|)$$
where $\alpha(t)=ae^R,$ $\displaystyle \gamma(t)=(1+\tan^2(1-c)R),$ and $\varphi(t)=(1-e^{-t})$ and $\psi(t)=\tan(1-c)t.$\\
Apply Theorem \ref{thm2},   (\ref{Ex4}), with $a$ small enough and $c=1-a,$   has a unique solution in $\Omega=\left\{x\in C([0,1]); \|x\|\leq 3\right\} $, in fact the solution is $x(t)=t.$

\noindent The following table show the numerical results of the proposed
 method for $a=0.01, b=1, R=3,$ and $x_0(t)=\sin(t).$
\begin{center}
	
	\small{\textsc{Table 4. Numerical results for    (\ref{Ex4})   with initial $x_0(t)=\sin(t).$}}

\begin{tabular}{|c|c|c|c|c|c|}
		\hline \hline
		& 	& \multicolumn{2}{c|}{$n_1=\dots=n_m=9 $ } & \multicolumn{2}{c|}{$n_1=\dots=n_m=33 $ }      \\
		\hline
		$t$  &$x^*(t)$	& $m=2$ & $m= 4$   & $m= 2$ & $m= 4$     \\
		\hline \hline
		$0.1$& $0.1$	&0.09994959265924196& 0.09994959275258121 & 0.09997341307990056&  0.09997341318295203 \\
		\hline
		$0.2$&  $0.2$	&  0.19982697772113361 &  0.19982698062053245&
 0.19994196511596454 &  0.19994196762406424  \\
		\hline
		$0.3$ &$0.3$	&0.29970145954436234 &  0.2997014781005956 &  0.2999105557764353  &  0.29991056948622924   \\
		\hline
		$0.4$& $0.4$	&
0.39957605185527684 &
0.3995761128223367 &
0.39987915932823637  &
0.3998792008487213  \\
		\hline
		$0.5$& $0.5$	& 0.49945067663057896  & 0.49945081633085925& 0.4998477565970089 &    0.49984784689621164  \\
		\hline
		$0.6$ & $0.6$	&  0.5993252823989378&  0.5993255387084228 &  0.5998163380147611 &    0.5998164954408373  \\
		\hline
		$0.7$& $0.7$	&
0.699199836738067  &
0.6992002390137386&
0.699784904433895&
0.6997851365136741  \\
		\hline
		$0.8$& $0.8$	& 0.799074326433192 & 0.7990748839377436& 0.7997534658105931 &   0.7997537620153589\\
		\hline
		$ 0.9$ &$ 0.9$	&  0.8989487530236073& 0.8989494465775325 &  0.8997220384049086 &  0.8997223654190059  \\
		\hline
		$ 1$ &$ 1$	&
0.9988231278772514   &
0.9988239054111422 &
0.9996906411587515 &
0.9996909415162489   \\
		
		\hline \hline
	\end{tabular}
\end{center}

\begin{center}
	\begin{tabular}{|c|c|c|c|c| }
		\hline \hline
		 	& \multicolumn{2}{c|}{$n_1=\dots=n_m=9 $ } & \multicolumn{2}{c|}{$n_1=\dots=n_m=33$ }  \\
		\hline
		  	& $m=2$ & $m= 4$   & $m= 2$ & $m= 4$    \\
		\hline \hline
		 $\|x^*-\tilde{x}\|_\infty$	& $1.17687\times 10^{-3}$ &  $1.17609\times 10^{-3}$  & $3.09359\times 10^{-4}$ & $3.09058\times 10^{-4}$  \\
		\hline \hline
	\end{tabular}
\end{center}
\end{exam}
\begin{exam}
 Consider the nonlinear differential equation
\begin{equation}\label{Ex5}
\displaystyle x(t)=  \displaystyle at\left[(b+t)^2+\frac{t}{(t+1)} \int_0^{t}\left(1-e^{-(t+1)(as+1)}\right)ds\right]^{-1} \left[(b+t)^2
+\int_0^t \int_0^{x(s)+1}e^{-(t+1)u}du ds\right].
\end{equation}
According to the above discussion,   (\ref{Ex5}) can be written as a fixed point problem $$x=F(x)\cdot G(x),$$ where $F$ and $G$ are defined in (\ref{EI}), with
$f(t,x)=\displaystyle at\left[(b+t)^2+\frac{t}{(t+1)} \int_0^{t}\left(1-e^{-(t+1)(as+1)}\right)ds\right]^{-1}$ and
$k(t,s,x)=\int_0^{x+1}e^{-(t+1)u}du.$

Let $0<R<1$ and let $x, y \in [-R,R].$  By an elementary calculus, we can show that
$$|f(t,x)-f(t,y)|\leq\alpha(t)\varphi(|x-y|)$$
and $$|k(t,s,x)-k(t,s,y)|\leq\gamma(t)\psi(|x-y|)$$
where  $ \alpha(t)=\varphi(t)=0,$ $ \psi(t)=\int_0^{2t}e^{-s}ds,$ and $\gamma(t,s)=\frac{1}{t+1}e^{(t+1)(R-1)}.$\\
 Taking $a=0.1, b=1,$ and applying Theorem \ref{thm2},  (\ref{Ex5}) has a unique solution  in $\Omega=\left\{x\in C([0,1]); \|x\|\leq R\right\}$. In fact the solution is $at.$
\begin{center}
	
 \small{\textsc{Table 5. Numerical results for    (\ref{Ex5})   with initial $x_0(t)=1/2 cos(10\pi t).$}}

	\begin{tabular}{|c|c|c|c|c|c| }
		\hline \hline
		& 	& \multicolumn{2}{c|}{$n_1=\dots=n_m=9 $ } & \multicolumn{2}{c|}{$n_1=\dots=n_m=33$ }  \\
		\hline
		$t$  &$x^*(t)$	& $m=2$ & $m= 4$   & $m= 2$ & $m= 4$    \\
		\hline \hline
		$0.1$& $0.01$	& 0.009807889768197995& 0.009807889768197995  & 0.009850176149181912 &  0.009850173620253975 \\
		\hline
		$0.2$&  $0.02$	&  0.01913347555848587 & 0.01913346934141619   &0.019763982715543943  &  0.01976400675926519  \\
		\hline
		$0.3$ &$0.03$	&  0.028858901595188682 &   0.028858870390823594& 0.029713698642337805 & 0.029713848529122386\\
		\hline
		$0.4$& $0.04$	&
0.038745648738524936 &
0.038745618536895766 & 0.039685125981635705  &
0.039685476825051164  \\
		\hline
		$0.5$& $0.05$	&  0.04868660113266646 &  0.0486866179763731  & 0.04967020617859288 &  0.04967087311797989 \\
		\hline
		$0.6$ & $0.06$	& 0.05866572236367756 &   0.05866579674631664  &0.05966446431697861  &  0.05966546941999516\\
		\hline
		$0.7$ & $0.06$	&
0.06866841551194598 &
0.06866853944486338 & 0.06966463823770887 &
0.06966603029961263    \\
		\hline
		$0.8$& $0.08$	&  0.07868630321311545&  0.07868650513410157 & 0.07966879607286238 &  0.07967053753105567 \\
		\hline
		$ 0.9$ &$ 0.09$	& 0.08871376330117915 &   0.08871405879246874 &0.08967554656783944  &  0.08967762811140002 \\
			\hline
		$1 $ &$ 0.09$	&
0.0987469779423797&
0.0987473453913395 & 0.09968401150389958 &
0.09968636366339986   \\
		\hline \hline
	\end{tabular}
\end{center}

\begin{center}
	\begin{tabular}{|c|c|c|c|c| }
		\hline \hline
		 	& \multicolumn{2}{c|}{$n_1=\dots=n_m=9 $ } & \multicolumn{2}{c|}{$n_1=\dots=n_m=33$ }  \\
		\hline
		  	& $m=2$ & $m= 4$   & $m= 2$ & $m= 4$    \\
		\hline \hline
		 $\|x^*-\tilde{x}\|_\infty$	& $1.33714\times 10^{-3}$ & $1.33705\times 10^{-3}$  &  $3.35272\times 10^{-4}$ &  $3.34982\times 10^{-4}$  \\
		\hline \hline
	\end{tabular}
\end{center}
\end{exam}

\section{\textbf{Conclusions}}{\label{sec:conclusions}

In this paper we have presented a  numerical method,  based } on the use of Schauder's bases, to solve hybrid nonlinear equations  in Banach algebras. To do this,  we have used  Boyd-Wong's  theorem to establish the existence and uniqueness of a fixed point for the product of two nonlinear  operators  in Banach algebra (Theorem \ref{eq}). The method is applied to a wide class of nonlinear integro-differential equations such as the ones we have illustrated  by means of several numerical examples.

The possibility of applying this process or a similar idea to other types of hybrid equations or systems of such equations is open and we hope to discuss this in the near future.

\section*{\textbf{Acknowledgements}}

The research of Aref Jeribi and Khaled Ben Amara has been partially supported by the University of Sfax.

\noindent The research of Maria Isabel Berenguer has been partially supported by Junta de Andalucia, Project FQM359 and by the ``Maria de Maeztu'' Excellence Unit IMAG, reference CEX2020-001105-M, funded by MCIN/AEI/10.13039/501100011033/ .

\end{document}